\documentclass[english,journal]{IEEEtran}
\usepackage[T1]{fontenc}
\usepackage[latin9]{inputenc}
\usepackage{verbatim}
\usepackage{amsmath}
\usepackage{amsthm}
\usepackage{amssymb}
\usepackage{xargs}[2008/03/08]

\makeatletter
  \theoremstyle{plain}
  \newtheorem{lem}{\protect\lemmaname}
  \theoremstyle{definition}
  \newtheorem{example}{\protect\examplename}
  \theoremstyle{plain}
  \newtheorem{thm}{\protect\theoremname}
  \theoremstyle{remark}
  \newtheorem{rem}{\protect\remarkname}
  \theoremstyle{plain}
  \newtheorem{cor}{\protect\corollaryname}

\usepackage{psfrag,mathptmx}
\usepackage{amssymb,latexsym,amsmath}
\usepackage{mathtools}
\usepackage{balance}
\usepackage{cite}
\usepackage{url}
\usepackage{dsfont}
\usepackage{bm}

\DeclareUrlCommand\email{\urlstyle{tt}\scriptsize}
\urldef{\raborges}{\email}{raborges@ene.unb.br}
\urldef{\ishihara}{\email}{ishihara@ene.unb.br}
\urldef{\kussaba}{\email}{kussaba@lara.unb.br}

\DeclareUrlCommand\site{\urlstyle{tt}\scriptsize}
\urldef{\YALMIP}{\site}{http://users.isy.liu.se/johanl/yalmip/}
\urldef{\SeDuMi}{\site}{http://sedumi.ie.lehigh.edu/}
\urldef{\PENLAB}{\site}{http://web.mat.bham.ac.uk/kocvara/penlab/}
\urldef{\ROLMIP}{\site}{http://www.dt.fee.unicamp.br/~agulhari/rolmip/rolmip.htm}

\makeatletter
\def\ps@headings{\def\@oddhead{\IEEEdoarxivheader{-1\oddsidemargin}\relax
\hbox{}\@IEEEheaderstyle\rightmark\hfil\thepage}\relax
\def\@evenhead{\IEEEdoarxivheader{-1\evensidemargin}\relax
\hbox{}\@IEEEheaderstyle\rightmark\hfil\thepage}\relax
\def\@oddfoot{\IEEEdoarxivfooter{-1\oddsidemargin}\hfil\hbox{}}\relax
\def\@evenfoot{\IEEEdoarxivfooter{-1\evensidemargin}\hfil\hbox{}}\relax}
\def\ps@IEEEtitlepagestyle{\ps@headings}

\def\IEEEarxivheadfootoffset{3pt}
\newdimen\IEEEheadtotopofpage
\newdimen\IEEEfoottobottomofpage
\newbox\@IEEEboxX

\def\IEEEarxivheader{}
\def\IEEEarxivfooter{}

\advance\IEEEheadtotopofpage by 1\topmargin
\advance\IEEEheadtotopofpage by 1in\relax

\IEEEfoottobottomofpage\paperheight
\advance\IEEEfoottobottomofpage by -1in\relax
\advance\IEEEfoottobottomofpage by -1\topmargin
\advance\IEEEfoottobottomofpage by -1\headheight
\advance\IEEEfoottobottomofpage by -1\headsep
\advance\IEEEfoottobottomofpage by -1\textheight
\advance\IEEEfoottobottomofpage by -1\footskip
\def\IEEEarxivheaderstyle{\normalfont\footnotesize}

\def\IEEEdoarxivheader#1{\@IEEEtrantmpdimenA\IEEEarxivheadfootoffset\relax
\@IEEEtrantmpdimenA -1\@IEEEtrantmpdimenA
\advance\@IEEEtrantmpdimenA by \IEEEheadtotopofpage
\settoheight{\@IEEEtrantmpdimenB}{\IEEEarxivheaderstyle HT}\relax
\advance\@IEEEtrantmpdimenA by -1\@IEEEtrantmpdimenB
\setbox\@IEEEboxX=\hbox{\relax
\raisebox{\@IEEEtrantmpdimenA}[0pt][0pt]{\parbox[t]{\textwidth}{\centering
\IEEEarxivheaderstyle\IEEEarxivheader}}}\relax
\wd\@IEEEboxX=0pt\relax
\ht\@IEEEboxX=0pt\relax
\dp\@IEEEboxX=0pt\relax
\box\@IEEEboxX\relax}

\def\IEEEarxivfooterstyle{\normalfont\footnotesize}

\def\IEEEdoarxivfooter#1{\@IEEEtrantmpdimenA\IEEEfoottobottomofpage\relax
\advance\@IEEEtrantmpdimenA by \IEEEarxivheadfootoffset\relax
\settodepth{\@IEEEtrantmpdimenB}{\IEEEarxivheaderstyle gjpqy}\relax
\advance\@IEEEtrantmpdimenA by 1\@IEEEtrantmpdimenB
\setbox\@IEEEboxX=\hbox{\hskip#1\hskip -1in\relax
\raisebox{\@IEEEtrantmpdimenA}[0pt][0pt]{\parbox[b]{\paperwidth}{\centering
\IEEEarxivfooterstyle\IEEEarxivfooter}}}\relax
\wd\@IEEEboxX=0pt\relax
\ht\@IEEEboxX=0pt\relax
\dp\@IEEEboxX=0pt\relax
\box\@IEEEboxX\relax}

\def\@IEEEheaderstyle{\normalfont\scriptsize}
\def\@IEEEfooterstyle{\normalfont\scriptsize}

\makeatother

\renewcommand{\IEEEarxivheadfootoffset}{3pt}

\renewcommand{\IEEEarxivheaderstyle}{\normalfont\footnotesize}
\renewcommand{\IEEEarxivfooterstyle}{\normalfont\footnotesize}

\renewcommand{\IEEEarxivheader}{This is the author's version of an article
that has been published in this journal. Changes were made to this version
by the publisher prior to publication.\\
The final version of record is available at 
\url{https://doi.org/10.1109/tac.2017.2682221}}
\renewcommand{\IEEEarxivfooter}{Copyright (c) 2017 IEEE. Personal use is
permitted. For any other purposes, permission must be obtained from the
IEEE by emailing pubs-permissions@ieee.org.}

\makeatother

\usepackage{babel}
\providecommand{\corollaryname}{Corollary}
\providecommand{\examplename}{Example}
\providecommand{\lemmaname}{Lemma}
\providecommand{\remarkname}{Remark}
\providecommand{\theoremname}{Theorem}

\begin{document}
\global\long\def\ker{\operatorname{Ker}}

\global\long\def\im{\operatorname{Im}}

\global\long\def\rank{\operatorname{rank}}

\global\long\def\argmin{\operatorname{argmin}}

\global\long\def\myvec#1{\ensuremath{#1}}

\global\long\def\mymatrix#1{\ensuremath{#1}}

\global\long\def\card{\operatorname{card}}

\newcommandx\simplex[1][usedefault, addprefix=\global, 1=N]{\ensuremath{\Delta_{#1}}}

\global\long\def\hppd#1{\ensuremath{\mymatrix{#1}(\myvec{\alpha})}}

\global\long\def\lpd#1{\ensuremath{#1(\myvec{\alpha})}}

\global\long\def\shpg#1{\ensuremath{\mathbb{H}_{\left(#1\right)}}}

\global\long\def\shppd#1#2#3{\ensuremath{\mathbb{H}_{\left(#1\right)}^{#2\times#3}}}

\global\long\def\ssym#1#2{\ensuremath{\mathbb{S}^{#1\times#2}}}

\global\long\def\bref#1{\mbox{(\ref{=00003D00003D00003D0000231})}}

\global\long\def\Hd{{\mathcal{H}}_{2}}

\global\long\def\Hi{{\mathcal{H}}_{\infty}}

\global\long\def\I{{\mbox{{\bf I}}}}

\global\long\def\Z{{\mbox{{\bf 0}}}}

\global\long\def\real{\mathds{R}}

\global\long\def\reais{\mathds{R}}

\global\long\def\C{\mathds{C}}

\global\long\def\N{\mathds{N}}

\global\long\def\inteiros{\mathbb{Z}}

\newcommandx\Zpos[1][usedefault, addprefix=\global, 1=]{\ensuremath{\mathbb{Z}_{+}^{#1}}}

\global\long\def\FinQForm{\text{(F1)}}

\global\long\def\FinExtScalar{\text{(F2)}}

\global\long\def\FinScalarP{\text{(F2a)}}

\global\long\def\FinScalarC{\text{(F2b)}}

\global\long\def\FinScalarRational{\text{(F2c)}}

\global\long\def\FinScalarPoly{\text{(F2d)}}

\global\long\def\FinScalarU{\text{(F2e)}}

\global\long\def\FinScalarf{\text{(F2f)}}

\global\long\def\FinScalarg{\text{(F2g)}}

\global\long\def\FinScalarh{\text{(F2h)}}

\global\long\def\FinExtMatrix{\text{(F3)}}

\global\long\def\FinMatrixP{\text{(F3a)}}

\global\long\def\FinMatrixC{\text{(F3b)}}

\global\long\def\FinMatrixU{\text{(F3c)}}

\global\long\def\FinOrthoComp{\text{(F4)}}

\title{Existence of continuous or constant Finsler's variables for parameter
dependent systems}

\author{J. Y. Ishihara, H. T. M. Kussaba, R. A. Borges\thanks{The authors are with the Department of Electrical Engineering, University
of Brasília, Brazil. A portion of this work was completed while J.
Y. Ishihara was a visiting associate researcher at Univ. California
Los Angeles (UCLA). \ishihara, \kussaba, \raborges.}}
\maketitle
\begin{abstract}
Finsler's lemma is a classic mathematical result with applications
in control and optimization. When the lemma is applied to parameter
dependent LMIs, as such those that arise from problems of robust stability,
the extra variables introduced by this lemma also become dependent
on this parameter. This technical note presents some sufficient conditions
which ensure, without losing generality, that these extra variables
can assume a simple functional dependence on the parameters as continuity
or even independence. The results allow avoiding an unnecessary use
of a more functionally complicated parameter dependent variable that
increases the search computational burden without reducing the conservatism
of the solution.
\end{abstract}


\maketitle

\section{Introduction }

Recently, the classical Finsler's lemma was reinterpreted to give
novel LMI characterizations to stability and control problems \cite{dOS:01}.
One important fact that motivates the search for LMI formulations
is that nowadays the LMI theory is equipped with a large number of
techniques which allows one to extend LMI analysis characterization
to control or filtering design conditions in a quite systematic way
(see, e.g. \cite{SIG:98}). The LMI theory has been successfully applied
to state a large number of problems arising from system and control
theory as standard convex optimization problems \cite{BEFB:94,EN:00,DY:13}.
Once a problem is stated as an LMI, it can be handled numerically
in polynomial time by several software such as the LMI Control Toolbox
\cite{GNLC:95} and SeDuMi \cite{Stu:99}, and can even, in some cases,
be solved analytically \cite{BEFB:94}.

Another advantage of the application of the Finsler's lemma is that
it enables one to deal with robustness issues using parameter dependent
(PD) formalism \cite{dOS:01}. It was observed that certain PD quadratic
or nonlinear matricial inequalities can be recast as PD-LMIs through
the Finsler's lemma. With this, several robust control problems turned
out to be analyzable in PD-LMI context without loss of generality.
Prior to the application of Finsler's lemma, in the literature there
were non-conservative but very computationally expensive nonlinear
inequality solutions or less expensive but approximated and conservative
LMI solutions. 

However, one should remember that although a PD-LMI formulation is
less computationally demanding than a nonlinear inequality formulation,
a PD-LMI is still very computationally demanding. It is known that
a generic PD-LMI is an NP-hard problem\cite{B-TN:98}. In face of
this, there have been several approaches in the linear matrix inequality
literature to relax PD-LMIs problems. For instance, one can consider
the worst case scenario for the effect of the parameter over a property
of interest, such as stability or $\mathcal{H}_{\infty}$ gain. Although
this approach is effective, it leads to very conservative criteria.
A less conservative but naive method to solve a PD-LMI would be to
approximate the parameter space $S$ by a finite set of points\textemdash this
would relax the PD-LMI into a finite number of LMIs. However, this
method usually requires a great number of points for its results to
become satisfactory, which generates a heavy computational load. Moreover,
if some points are excluded, it is possible that the result turns
out to be excessively optimist \cite{AT:00}. Another relaxation procedure
is to use a polynomial relaxation, which consists of restricting the
variables of the PD-LMI to polynomial functions of a fixed degree
$g$ and using the matrix coefficients of the polynomial stemmed from
this procedure to generate new LMIs independent of the parameter.
 With the increase of the degree $g$ of the polynomial, less conservative
sets of conditions are obtained. From degrees greater than some value,
the set of conditions imply the original PD-LMI \cite{Bli:04b}.

Along with the technique of polynomial relaxation, in some PD-LMI
solvable problems, one may also use the Finsler's lemma to relax conservativeness
in further analysis. For instance, in the context of robust filter
or control design, consider an uncertain continuous-time linear system
given by 
\begin{equation}
\dot{x}(t)=A(s)x(t)\label{eq:uncertain_system}
\end{equation}
where ${s}$ is a vector of unknown parameters belonging to a known
set $S$. It is possible to prove (see \cite{dOS:01}) that a sufficient
condition to robust stability is related to the existence of a positive
definite matrix valued function $P(s)$ and a scalar function $\mu(s)$
such that for all $s\in S$, 
\begin{equation}
\begin{bmatrix}-\mu(s)A^{T}(s)A(s) & \mu(s)A^{T}(s)+P(s)\\
\mu(s)A(s)+P(s) & -\mu(s)I
\end{bmatrix}\prec0.\label{eq:stability_scalar_Finsler_mu_function}
\end{equation}
Differently from the classical approach, in (\ref{eq:stability_scalar_Finsler_mu_function}),
the matrix variable $P(s)$ does not multiply the system matrix $A(s)$.
This separation allows one to deal with stability analysis of closed
loop systems to design stabilizing controllers or filters in a much
simpler procedure than the traditional. Also, it  may lead to less
conservative results \cite{OP:07a,EPA:15}. As a consequence, Finsler's
lemma has been notably employed in the literature in several contexts\cite{EHPA:05,SFS:06,MOP:09a,AS:11,ELG:14,FLOP:15,GE-ML:15}.
For examples of several control problems that can be dealt with the
Finsler's lemma (also known as S-lemma) one can see \cite{EPA:15}.
The drawback of this approach is that it increases the search space.
For instance, in contrast to the traditional stability analysis with
a unique variable $P(s)$, (\ref{eq:stability_scalar_Finsler_mu_function})
is a PD-LMI in $P(s)$ and $\mu(s)$, that is, now the introduced
extra function $\mu(s)$ is also a variable that needs to be found.

Relaxation is a good resort when the solution space is little known.
Another approach considered in the literature is to find some properties
on the parameter set $S$ and on the matrices of the system in order
to reduce the search space of the solutions \cite{Bli:04b,Cim:15}.
In fact, if $S$ is a compact set and the matrices of the PD-LMI depend
continuously on the parameter, in \cite{Bli:04b} it is proved that
if the PD-LMI has a solution, then one can restrict the search for
a solution in the set of polynomial functions. In the particular case
that $S$ is a unit simplex, there will also be without loss of generality,
a solution that is a homogeneous polynomial \cite{BOMP:06b}. These
results are important since, with the knowledge that a continuous
function can be approximated by a polynomial, they form the theoretical
justification for the polynomial relaxation approach (where at least
a continuous solution is supposed to exist).  If the parameter space
is instead all the space and the matrices of the PD-LMI depends polynomially
on the parameter, then \cite{Cim:15} shows that one can restrict
the search space and look for a rational solution. 

Besides the polynomial or rational structures for the slack variable
introduced by the Finsler's lemma, we show in this paper that there
are some cases in which a parameter independent slack variable is
as good as a parameter dependent one. In other words, there is no
gain to impose a complicated functional structure to the extra variable.
In fact, the objective of this paper is to investigate some conditions
in which it is possible, without loss of generality, to impose solutions
of simpler functional structures, such as continuous, rational or
constant solutions for the slack variables introduced by Finsler's
lemma. In the LMI context, this allows reducing computational burden
without increasing the conservatism of the solution. This work is
an extension of our previous analysis on uniform versions of Finsler's
lemma \cite{KIB:15}. Here, some new results are presented and some
results from \cite{KIB:15} were improved: the converse of Lemma~3
from \cite{KIB:15} has been proved in Lemma~\ref{lem:sup_inf_uniform}
of this paper. Both Lemma~3 and Theorem~2 from \cite{KIB:15} were
stated for matrix-valued functions of a compact subset of $\mathbb{R}$.
In this paper, these results are extended for matrix-valued functions
of a compact subset of $\mathbb{R}^{d}$. The hypothesis of Theorem~1
and Corollary~1 from \cite{KIB:15} has been weakened in Theorem~\ref{thm:cont_Q_B}
of this work. 

\textit{Notation}: In the sequel the following notation will be used:
$\mathbb{R}$ is the set of real numbers, $\mathbb{R}^{m\times n}$
the set of real matrices of order $m\times n$. $\mathbb{S}^{n}$,
$\mathbb{S}_{>0}^{n}$ and $\mathbb{S}_{\ge0}^{n}$ are, respectively,
the set of symmetric, positive definite and positive-semidefinite
matrices of order $n\times n$. For a matrix $A$, $A\prec0$ indicates
that $-A\in\mathbb{S}_{>0}^{n}$, $A^{T}$ its transpose; $\im(A)$
and $\ker(A)$ are respectively the image and the kernel of $A$;
$A^{\perp}$ is a matrix whose columns span a basis for $\ker(A)$
and $A^{1/2}$ denotes the principal square root of a positive-semidefinite
matrix $A$. 

\section{The Finsler's lemma\label{sec:Main_results}}

The importance of the Finsler's lemma can be highlighted by the fact
that it is equivalent to other important results in control and optimization
literature such as Yakubovich's S-lemma \cite{Z-ZJ-H:10}. In fact,
Finsler's lemma have been proved several times \cite{Alb:38,Rei:38,AGH:73,Uhl:79,Ham:99,dOS:01,PT:07}.
Nowadays the Finsler's lemma is commonly stated as below:
\begin{lem}
\label{lem:Finsler_from_dOS:01} \cite{dOS:01} Let ${Q\in\mathbb{S}^{n}}$
and ${B\in\mathbb{R}^{m\times n}}$, with $\rank(B)<n$. Then the
following statements are equivalent: 

\begin{enumerate}
\item \label{enu:FinslerQuadForm} $x^{T}Qx<0$ for all $x\in\mathbb{R}^{n}$
such that $x\ne0$ and $Bx=0$. 
\item \label{enu:FinslerExtraScalar} There exists $\mu\in\mathbb{R}$ such
that $Q-\mu B^{T}B\prec0$. 
\item \label{enu:FinslerExtraMatrix} There exists $X\in\mathbb{R}^{n\times m}$
such that $Q+XB+B^{T}X^{T}\prec0.$
\item \label{enu:FinslerOrthoComp} $\left(B^{\perp}\right)^{T}QB^{\perp}\prec0.$ 
\end{enumerate}
\end{lem}
The equivalence between \ref{enu:FinslerQuadForm}) and \ref{enu:FinslerExtraScalar})
is attributed to Paul Finsler \cite{Fin:37}, where he considered
a more general case involving an indefinite matrix instead of the
positive semi-definite matrix $B^{T}B$. It is also interesting to
remark that the equivalence between \ref{enu:FinslerExtraMatrix})
and \ref{enu:FinslerOrthoComp}) can be seen as a particular case
of the Projection Lemma (also known as Elimination Lemma), which is
also widely used in control \cite{BEFB:94,PDSV:09}. 

The first contribution of this paper is the observation that in statements
\ref{enu:FinslerExtraScalar}) and \ref{enu:FinslerExtraMatrix})
of Lemma~\ref{lem:Finsler_from_dOS:01}, the hypothesis $\rank(B)<n$
is not necessary. 
\begin{lem}
\label{lem:Finsler_full_rank} Let ${Q\in\mathbb{S}^{n}}$ and ${B\in\mathbb{R}^{m\times n}}$,
with $\rank(B)=n$. Then there exists $\mu\in\mathbb{R}$ such that
$Q-\mu B^{T}B\prec0$. One such $\mu$ is given by
\[
\mu=1+\frac{\lambda_{max}(Q)+|\lambda_{max}(Q)|}{\lambda_{min}(B^{T}B)}.
\]

Moreover, there exists $X\in\mathbb{R}^{n\times m}$ such that $Q+XB+B^{T}X^{T}\prec0$.
One such $X$ is given by
\[
X=-\frac{1}{2}\left[1+\frac{\lambda_{max}(Q)+|\lambda_{max}(Q)|}{\lambda_{min}(B^{T}B)}\right]B^{T}.
\]
\end{lem}
\begin{IEEEproof}
See Lemma~\ref{lem:sufficient_condition_for_mu_continuous} in Appendix
and consider $M=Q$ and $N=B^{T}B$.
\end{IEEEproof}
With Lemma~\ref{lem:Finsler_full_rank}, one can generalize the Finsler's
lemma by eliminating the hypothesis $\rank(B)<n$. When $\rank(B)=n$,
the empty sentences \ref{enu:FinslerQuadForm}) and \ref{enu:FinslerOrthoComp})
are reinterpreted as being trivially satisfied. This general Finsler's
Lemma, extended to parameter dependent systems, is stated in the next
section. 

\section{Finsler's lemmas for parameter dependent systems }

When dealing with uncertain systems, the matrices $Q$ or $B$ can
become dependent on parameters \cite{OP:07a}. A pointwise extension
of Finsler's lemma is easily obtained as stated in the following. 
\begin{lem}
\label{lem:Finsler_pointwise} Let $S\subseteq\mathbb{R}^{d}$, ${Q:S\rightarrow\mathbb{S}^{n}}$
and ${B:S\rightarrow\mathbb{R}^{m\times n}}$. Then the following
statements are equivalent: \renewcommand{\theenumi}{\roman{enumi}}

\begin{enumerate}
\item[\textbf{$\FinQForm$}]  For each $s\in S$, one has $x^{T}Q\left(s\right)x<0$ for all $x\in\mathbb{R}^{n}$
such that $x\ne0$ and $B(s)x=0$. 
\item[\textbf{$\FinExtScalar$}] $\left(\forall s\in S\right)\left(\exists\mu\left(s\right)\in\mathbb{R}\right):Q\left(s\right)-\mu\left(s\right)B^{T}\left(s\right)B\left(s\right)\prec0$. 
\item[\textbf{$\FinExtMatrix$}]  For each $s\in S$, there exists $X(s)\in\mathbb{R}^{n\times m}$
such that
\[
Q(s)+X(s)B(s)+B^{T}(s)X^{T}(s)\prec0.
\]
\item[\textbf{$\FinOrthoComp$}]  For each $s\in S$, one has that $\left(B^{\perp}(s)\right)^{T}Q(s)B^{\perp}(s)\prec0.$ 
\end{enumerate}
\end{lem}
\begin{IEEEproof}
Follows directly extending pointwisely Lemma~\ref{lem:Finsler_from_dOS:01}
and Lemma~\ref{lem:Finsler_full_rank}.
\end{IEEEproof}
The inequalities in sentences $\FinExtScalar$ and $\FinExtMatrix$
are known as parameter dependent LMIs (PD-LMIs), and must be satisfied
for all parameters $s\in S$. Since the problem of finding a solution
to a PD-LMI may be NP-hard {[}10{]}, it is interesting to reduce the
search space by seeking solutions in classes of functions $\mu\left(s\right)$
or $X\left(s\right)$ with some functional structure like continuity,
rational or polynomial dependency or even independent constant solution.
Since we can obtain the solution $X\left(s\right)=-\frac{1}{2}\mu\left(s\right)B^{T}\left(s\right)$
for $\FinExtMatrix$ once we have a solution $\mu\left(s\right)$
for $\FinExtScalar$, in this paper we will be mainly interested in
the following situations for $\FinExtScalar$: 
\begin{enumerate}
\item[\textbf{$\FinScalarP$}]  There exists a function $\mu:S\rightarrow\mathbb{R}$ such that
$Q(s)-\mu(s)B^{T}(s)B(s)\prec0$ for all $s\in S$.
\item[\textbf{$\FinScalarC$}]  There exists a continuous function $\mu:S\rightarrow\mathbb{R}$
such that $Q(s)-\mu(s)B^{T}(s)B(s)\prec0$ for all $s\in S$.
\item[\textbf{$\FinScalarRational$}]  There exists a rational function $\mu:S\rightarrow\mathbb{R}$ without
singularities on $S$ such that $Q(s)-\mu(s)B^{T}(s)B(s)\prec0$ for
all $s\in S$.
\item[\textbf{$\FinScalarPoly$}]  There exists a polynomial function $\mu:S\rightarrow\mathbb{R}$
such that $Q(s)-\mu(s)B^{T}(s)B(s)\prec0$ for all $s\in S$.
\item[\textbf{$\FinScalarU$}]  There exists a constant $\bar{\mu}\in\mathbb{R}$ such that $Q\left(s\right)-\bar{\mu}B^{T}\left(s\right)B\left(s\right)\prec0$
for all $s\in S$.
\end{enumerate}
It is obvious that $\FinScalarU\Rightarrow\FinScalarPoly\Rightarrow\FinScalarRational\Rightarrow\FinScalarC\Rightarrow\FinScalarP$,
however, as shown in the following counterexample, reverse implications
are not true in general and additional hypothesis on the set $S$
and the functions $Q\left(s\right)$ and $B\left(s\right)$ are necessary
to assure that existence of pointwise solutions $\mu\left(s\right)$
guarantees the existence of a solution with a simple functional dependence
on $s$.
\begin{example}
\label{example:not_uniform_in_general-1} For $S=(0,+\infty)\subset\mathbb{R}$,
let $Q:S\rightarrow\mathbb{R}^{2\times2}$ and $B:S\rightarrow\mathbb{R}^{1\times2}$
be given by
\[
Q\left(s\right)=\begin{bmatrix}-1 & 0\\
0 & q\left(s\right)
\end{bmatrix},\ B^{T}\left(s\right)=\begin{bmatrix}0\\
s
\end{bmatrix}.
\]
One has that any solution $\mu\left(s\right)$ for $Q\left(s\right)-\mu\left(s\right)B^{T}\left(s\right)B\left(s\right)\prec0$
must satisfy
\begin{equation}
\mu\left(s\right)>\frac{q\left(s\right)}{s^{2}},\:s\in(0,+\infty).\label{eq:aux}
\end{equation}
Thus, independently if the function $q\left(s\right)$ is discontinuous
or not, to guarantee that $\FinScalarP$ will be satisfied one can
take $\mu\left(s\right)=\frac{q\left(s\right)}{s^{2}}+\epsilon$ with
$\epsilon>0$. 

For $q\left(s\right)=s$, one function satisfying $\FinScalarC$ and
$\FinScalarRational$ is $\mu\left(s\right)=\frac{1}{s}+\epsilon$
with $\epsilon>0$. It is easy to see that there is no constant solution
$\mu\left(s\right)=\bar{\mu}$ nor a polynomial solution $\mu\left(s\right)=a_{n}s^{n}+...+a_{1}s+a_{0}$
since in (\ref{eq:aux}) the function $\frac{q\left(s\right)}{s^{2}}=\frac{1}{s}$
grows without bound as $s$ goes to $0$. 

Similarly, for $q\left(s\right)=e^{s}$, one function satisfying $\FinScalarC$
is $\mu\left(s\right)=\frac{e^{s}}{s^{2}}+\epsilon$ with $\epsilon>0$.
And it can be shown that there is no constant, polynomial nor rational
solution $\mu\left(s\right)$, since otherwise from (\ref{eq:aux})
there would exist a rational function growing faster than the exponential
as $s\rightarrow\infty$.

Further, take a function $q(s)$ such that, for some point $0<\bar{s}<+\infty$,
$\frac{q\left(s\right)}{s^{2}}$ goes to $+\infty$ for $s\downarrow\bar{s}$
and goes to any value $\bar{q}<+\infty$ for $s\uparrow\bar{s}.$
In this case, it is impossible to find a continuous $\mu\left(s\right)$
satisfying the inequality (\ref{eq:aux}), that is, although $\FinScalarP$
is satisfied, $\FinScalarC$ is not. \qed
\end{example}
The main objective of this paper is to investigate under which conditions
there are equivalences among the statements $\FinScalarP$,$\ldots$,$\FinScalarU$
since this means that one can reduce the search space to subspaces
of functions $\mu(s)$ with a simpler structure without loss of generality.
Some of these equivalences have already been shown in the literature. 
\begin{thm}
\label{thm:Finsler_lemma_for_matrix_polynomials} \cite{Cim:15} If
the functions $Q\left(s\right)$ and $B\left(s\right)$ are polynomial
(matrices) over $S=\mathbb{R}^{d}$, then $\FinScalarP$ is equivalent
to $\FinScalarPoly$.
\end{thm}
\begin{IEEEproof}
Follows directly by choosing $F=-Q$ and $G=-B^{T}B$ in Proposition~3.2
from \cite{Cim:15}\footnote{The original Proposition 3.2 establishes the existence of $\mu$ rational
without singularities. However, from its proof, one can conclude that
$\mu$ can be chosen polynomial. We thank the anonymous reviewer that
pointed out the possibility of this extension and to prof. Jaka Cimpri\v{c}
for ratifying this fact.}.
\end{IEEEproof}
\begin{thm}
\label{thm:Bliman} \cite{Bli:04b} Let ${Q:S\rightarrow\mathbb{S}^{n}}$
and ${B:S\rightarrow\mathbb{R}^{m\times n}}$ be continuous matrix
valued functions on a compact $S\subset\mathbb{R}^{d}$. Then $\FinScalarP$
is equivalent to $\FinScalarPoly$.
\end{thm}
\begin{IEEEproof}
Follows directly from Theorem~1 of \cite{Bli:04b}.
\end{IEEEproof}
In the next section, we investigate other conditions for some equivalences
among the statements $\FinScalarP$,$\ldots$,$\FinScalarU$. In particular,
one of the main goals of this paper is to investigate when the PD-LMI
in (F2) is also valid uniformly in $\mu$, that is, when $\FinScalarP$
is equivalent to $\FinScalarU$.

\subsection{\label{subsec:results_for_scalar_mu} Results on the scalar-valued
function $\mu\left(s\right)$}

The first result of this paper shows that if the parameter set $S$
is compact, then there is no gain in searching for a complicated continuous
solution $\mu(s)$ since one can reduce the search to constant solutions.
\begin{thm}
\label{thm:Finsler_lemma_for_continuous_mu} Let ${Q:S\rightarrow\mathbb{S}^{n}}$
and ${B:S\rightarrow\mathbb{R}^{m\times n}}$ be matrix valued functions
on a compact $S\subset\mathbb{R}^{d}$. Then the sentences from $\FinScalarC$
to $\FinScalarU$ are all equivalent. If further, $Q$ and $B$ are
continuous then the sentences from $\FinScalarP$ to $\FinScalarU$
are all equivalent.
\end{thm}
\begin{IEEEproof}
The proof for $\FinScalarU\Rightarrow\FinScalarPoly\Rightarrow\FinScalarRational\Rightarrow\FinScalarC$
is immediate. Suppose now that $\FinScalarC$ is valid. By Weierstrass'
theorem \cite[p.90]{Sund:96}, the function $\mu\left(s\right)$ has
a maximum $\mu_{max}$ in $S$. Thus
\[
0\succ Q(s)-\mu(s)B^{T}(s)B(s)\succeq Q(s)-\mu_{max}B^{T}(s)B(s)
\]
 and $\FinScalarU$ is satisfied with $\bar{\mu}=\mu_{max}$.  If
$Q$ and $B$ are continuous, then $\FinScalarP$ is equivalent to
$\FinScalarPoly$ by Theorem~\ref{thm:Bliman}. Since all sentences
from $\FinScalarC$ to $\FinScalarU$ are equivalent, the result follows.
\end{IEEEproof}
Note that for PD-LMIs with one variable $\mu\left(s\right)$, Theorem
\ref{thm:Finsler_lemma_for_continuous_mu} extends the consecrated
result of \cite{Bli:04b}. This might be useful in the context of
polynomial relaxation procedures for PD-LMIs \cite{Bli:04b,BOMP:06b},
wherein it is assumed that $\mu\left(s\right)$ is a polynomial function.
With Theorem \ref{thm:Finsler_lemma_for_continuous_mu}, one can reduce
the search space from the set of polynomial functions to the set of
real numbers. 

Another fundamental lemma which gives a general condition to guarantee
the existence of one $\mu$ for all $s\in S$ is presented below.
\begin{lem}
\label{lem:sup_inf_uniform}Let ${Q:S\rightarrow\mathbb{S}^{n}}$
and ${B:S\rightarrow\mathbb{R}^{m\times n}}$ be functions on $S\subset\mathbb{R}^{d}$
such that
\begin{equation}
\sup_{s\in S}\,\inf\left\{ \mu\in\mathbb{R}\mid Q\left(s\right)-\mu B^{T}\left(s\right)B\left(s\right)\prec0\right\} <\infty.\label{eq:lem_sup-inf_sup-inf}
\end{equation}

Then the statements from $\FinScalarP$ to $\FinScalarU$ are all
equivalent. Conversely, if $\FinScalarU$ holds, then $\FinScalarP$
and (\ref{eq:lem_sup-inf_sup-inf}) hold.
\end{lem}
\begin{IEEEproof}
The proof for $\FinScalarU\Rightarrow\FinScalarPoly\Rightarrow\FinScalarRational\Rightarrow\FinScalarC\Rightarrow\FinScalarP$
is immediate. We now prove that $\FinScalarP$ implies $\FinScalarU$.
For each $s\in S$, define
\[
\mathcal{M}\left(s\right)=\left\{ \mu\in\mathbb{R}\mid Q\left(s\right)-\mu B^{T}\left(s\right)B\left(s\right)\prec0\right\} .
\]

For each $s$, one has that if $\mu^{*}\in\mathcal{M}\left(s\right)$,
then ${\mu^{*}+\alpha\in\mathcal{M}\left(s\right)}$ for all $\alpha\ge0$,
since
\[
0\succ Q\left(s\right)-\mu^{*}B^{T}\left(s\right)B\left(s\right)\succeq Q\left(s\right)-\left(\mu^{*}+\alpha\right)B^{T}\left(s\right)B\left(s\right).
\]

By $\FinScalarP$, one has that $\mathcal{M}\left(s\right)\ne\emptyset$
for any $s\in S$. Therefore, there always exists a $\mu^{*}\left(s\right)\in\mathcal{M}\left(s\right)$
such that $\left[\mu^{*}\left(s\right),+\infty\right)\subseteq\mathcal{M}\left(s\right)$.
By (\ref{eq:lem_sup-inf_sup-inf}), there exists $m\in\mathbb{R}$
such that
\[
\inf\left\{ \mu\in\mathbb{R}\mid Q\left(s\right)-\mu B^{T}\left(s\right)B\left(s\right)\prec0\right\} <m,\ \forall s\in S.
\]

Consequently, $m\in\mathcal{M}\left(s\right)$ for all $s\in S$ and
$\bar{\mu}=m$ is such that $\FinScalarU$ holds. Conversely, suppose
that $\FinScalarU$ holds. It is clear that $\FinScalarP$ holds and
$\mathcal{M}\left(s\right)\ne\emptyset$ for all $s\in S$ since $\bar{\mu}\in\mathcal{M}\left(s\right)$.
It follows that, $\inf\mathcal{M}\left(s\right)\leq\bar{\mu},\ \forall s\in S$
and $\sup_{s\in S}\,\inf\mathcal{M}\left(s\right)\leq\bar{\mu}<\infty.$
\end{IEEEproof}
\begin{rem}
If one adopts the convention that $\inf\emptyset=+\infty$ then Lemma~\ref{lem:sup_inf_uniform}
can be restated as: $\FinScalarU$ is equivalent to (\ref{eq:lem_sup-inf_sup-inf}).
The above presentation of the lemma was preferred in order to distinguish
the situations where does not exist $\mu\left(s\right)$ from those
where the function $\mu\left(s\right)$ goes to infinity.
\end{rem}
Note that Lemma~\ref{lem:sup_inf_uniform} is very general and does
not require any special structure on the functions $Q\left(s\right)$
and $B\left(s\right)$ nor on the set $S$. One immediate use of Lemma~\ref{lem:sup_inf_uniform}
is a derivation of a Finsler's lemma version for switching systems.
This is a case where (\ref{eq:lem_sup-inf_sup-inf}) can be easily
checked as discussed in the following.
\begin{cor}
\label{cor:swiching}Consider matrix valued functions $Q:S\rightarrow\mathbb{S}^{n}$
and $B:S\rightarrow\mathbb{R}^{m\times n}$ assuming a finite number
of values given by the set $\left\{ \left(Q_{1},B_{1}\right),...,\left(Q_{N},B_{N}\right)\right\} $.
Then the statements from $\FinScalarP$ to $\FinScalarU$ and the
following statements are equivalent:

\begin{enumerate}
\item[$\FinScalarf$]  There exists a constant $\bar{\mu}\in\mathbb{R}$ such that $Q_{i}-\bar{\mu}B_{i}^{T}B_{i}\prec0$
for every $i\in\left\{ 1,...,N\right\} $.
\item[$\FinScalarg$]  For every $i\in\left\{ 1,...,N\right\} ,\exists\mu_{i}\in\mathbb{R}$
such that $Q_{i}-\mu_{i}B_{i}^{T}B_{i}\prec0$.
\end{enumerate}
Furthermore, if $Q\left(s\right)$ and $B\left(s\right)$ are piecewise
constant functions on $S$, then $\left(F2a\right)-\left(F2g\right)$
are also equivalent to

\begin{enumerate}
\item[$\FinScalarh$]  There exists a piecewise constant function $\mu:S\rightarrow\mathbb{R}$
such that $Q\left(s\right)-\mu\left(s\right)B^{T}\left(s\right)B\left(s\right)\prec0$
for all $s\in S$.
\end{enumerate}
\end{cor}
\begin{IEEEproof}
Consider $\FinScalarP$ valid. It is clear that for every $i\in\left\{ 1,...,N\right\} $,
$\mathcal{M}_{i}\coloneqq\left\{ \mu\in\mathbb{R\mid}Q_{i}-\mu B_{i}^{T}B_{i}\prec0\right\} \neq\emptyset$
and so $\bar{\mu}_{i}\coloneqq\inf\mathcal{M}_{i}<+\infty$. Then,
$\sup_{s\in S}\,\inf\left\{ \mu\in\mathbb{R}\mid Q\left(s\right)-\mu B^{T}\left(s\right)B\left(s\right)\prec0\right\} =\sup\left\{ \bar{\mu}_{1},...,\bar{\mu}_{N}\right\} <+\infty.$
From Lemma~\ref{lem:sup_inf_uniform} it follows that $\FinScalarP\implies\FinScalarU$.
The implications $\FinScalarU\Rightarrow\FinScalarf\Rightarrow\FinScalarg\Rightarrow\FinScalarP$
and $\FinScalarPoly\Rightarrow\FinScalarRational\Rightarrow\FinScalarC\Rightarrow\FinScalarP$
are immediate. For the last statement of the corollary, define the
set $S_{i}\coloneqq\left\{ s\in S\mid\left(Q\left(s\right),B\left(s\right)\right)=\left(Q_{i},B_{i}\right)\right\} $,
consider any $\mu_{i}\in\left\{ \mu\in\mathbb{R\mid}Q_{i}-\mu B_{i}^{T}B_{i}\prec0\right\} $
and $\mu\left(s\right)=\mu_{i},\,s\in S_{i}$.
\end{IEEEproof}
One particular case of interest for the Corollary~\ref{cor:swiching}
is when the system parameters $Q\left(s\right)$ and $B\left(s\right)$
represent a switching system, and the parameter set is interpreted
as continuous time, $S=\mathbb{R}$, or as discrete time, $S=\mathbb{Z}$.
Note that if $\left(Q_{i},B_{i}\right)$ are appropriately chosen,
then the Corollary~\ref{cor:swiching} states that, no matter which
switching policy is used, there always exists a continuous or even
a constant solution $\mu$. The derivation of Finsler's lemma versions
for switching systems is now immediate. One such lemma is given in
the following.
\begin{lem}
Consider matrix sequences $Q\left(k\right)\in\left\{ Q_{1},...,Q_{N_{Q}}\right\} \subset\mathbb{S}^{n}$
and $B\left(k\right)\in\left\{ B_{1},...,B_{N_{B}}\right\} \subset\mathbb{R}^{m\times n}$,
$k=1,2,...$. Then the following statements are equivalent:

\begin{enumerate}
\item[$\left(i\right)$]  For each $k=1,2,...$, one has $x^{T}Q\left(k\right)x<0$ for all
$x\in\mathbb{R}^{n}$ such that $x\ne0$ and $B(k)x=0$.
\item[$\left(ii\right)$]  There exists a constant $\bar{\mu}\in\mathbb{R}$ such that $Q_{i}-\bar{\mu}B_{j}^{T}B_{j}\prec0$
for every $i\in\left\{ 1,...,N_{Q}\right\} $ and for every $j\in\left\{ 1,...,N_{B}\right\} $.
\end{enumerate}
\end{lem}
The next theorem gives conditions for the existence of a continuous
solution when it is known that a pointwise solution exists.
\begin{thm}
\label{thm:cont_Q_B} Let ${Q:S\rightarrow\mathbb{S}^{n}}$ and ${B:S\rightarrow\mathbb{R}^{m\times n}}$
be continuous functions on $S\subseteq\mathbb{R}^{d}$. Then the statements
$\FinScalarP$ and $\FinScalarC$ are equivalent. If further, $S$
is a compact set, then the statements from $\FinScalarP$ to $\FinScalarU$
are all equivalent.
\end{thm}
\begin{IEEEproof}
Define the function $\mu_{\inf}:S\rightarrow\mathbb{R}\cup\left\{ -\infty\right\} $
as $\mu_{\inf}(s)\coloneqq\inf\,\mathcal{M}\left(s\right)$ where
$\mathcal{M}\left(s\right)\coloneqq\left\{ \mu\in\mathbb{R}\mid Q\left(s\right)-\mu B^{T}\left(s\right)B\left(s\right)\prec0\right\} $.
By (F2a), that is, non emptiness of $\mathcal{M}\left(s\right),$
it is clear that $\mu_{\inf}\left(s\right)<+\infty$ for all $s\in S$
and $\mu_{\inf}(s)=-\infty$ for $s\in S$ at which $Q\left(s\right)\prec0$
and $B\left(s\right)=0.$ Since $\bar{\mu}\in\mathcal{M}\left(s\right)$
implies that $\bar{\mu}+\alpha\in\mathcal{M}\left(s\right)$ for any
$\alpha>0$, the set $\mathcal{M}\left(s\right)$ is an open interval
of the form $\mathcal{M}\left(s\right)=\left(\mu_{\inf}(s),+\infty\right).$
Taking the induced topology of $\mathbb{R}^{d}$ on $S$, one can
use the same arguments in the proof for continuity on $\mathbb{R}^{d}$
of Lemma~3.1 of \cite{Cim:15} to conclude that the extended real-valued
function $\mu_{\inf}$ is continuous on $S$. 

Now, take any $\varepsilon>0$ and define the real valued function
$\tilde{\mu}(s)\coloneqq\max\left\{ \mu_{\inf}(s)+\varepsilon,0\right\} $.
It is clear that $\tilde{\mu}:S\rightarrow\mathbb{R}$ is a continuous
function with $\tilde{\mu}(s)\in\mathcal{M}\left(s\right)$ for all
$s\in S$. Thus (F2b) is satisfied with $\tilde{\mu}$.
\end{IEEEproof}
Theorem \ref{thm:cont_Q_B} states that a continuous solution exists
if $Q$ and $B$ are continuous over an arbitrary subset $S\subseteq\mathbb{R}^{d}$.
In contrast, Theorem~\ref{thm:Bliman} allows a polynomial solution
but just when $S$ is compact. Theorem~\ref{thm:Finsler_lemma_for_matrix_polynomials}
allows a polynomial solution if the functions $Q$ and $B$ are polynomials
over $S=\mathbb{R}^{d}$.

The next example illustrates an application of Theorem~\ref{thm:cont_Q_B}
in the context of non-linear systems.
\begin{example}
\label{exa:-ExponentialStability}Consider
\begin{equation}
\dot{x}=f(x,t)+B^{T}(x,t)u,\:t>0\label{eq:non_linear_system}
\end{equation}
where $x\in\mathbb{R}^{n}$ and $u\in\mathbb{R}^{m}$ are the state
and the control variables, respectively. The functions $f:\mathbb{R}^{n}\times\mathbb{R}^{+}\rightarrow\mathbb{R}^{n}$
and $B:\mathbb{R}^{n}\times\mathbb{R}^{+}\rightarrow\mathbb{R}^{m\times n}$
are assumed to be smooth. In \cite{MTS:15} it is shown that system
(\ref{eq:non_linear_system}) is universally exponentially stabilizable
with rate $\lambda$ if there exists a positive definite matrix valued
function $M:\mathbb{R}^{n}\times\mathbb{R}^{+}\rightarrow\mathbb{S}_{>0}^{n}$
and a function $\rho:\mathbb{R}^{n}\times\mathbb{R}^{+}\rightarrow\mathbb{R}$
such that
\begin{equation}
\dot{M}+\left(\frac{\partial f}{\partial x}\right)^{T}M+M\left(\frac{\partial f}{\partial x}\right)-\rho MB^{T}BM+2\lambda M\preceq0\label{eq:DLMI_condition}
\end{equation}
for all $x\in\mathbb{R}^{n}$ and $t>0.$ 

Since the direct solution to  (\ref{eq:DLMI_condition}) is very hard
to be obtained, for simplicity, one may consider a restricted search
space by choosing $M=I$. This results in
\begin{equation}
Q(x)-\rho(x)B^{T}(x)B(x)\preceq-2\lambda I\prec0\label{eq:LMI_from_DMI}
\end{equation}
with
\[
Q(x)\coloneqq\left(\frac{\partial f}{\partial x}\right)^{T}+\left(\frac{\partial f}{\partial x}\right).
\]
Note that even in this simplified context, the search space is still
very big since, in principle, $\rho(x)$ could have any structure.
Theorem~\ref{thm:cont_Q_B} can be called to reduce the search space
since it guarantees that there is no loss of generality to assume
that a continuous solution $\rho(x)$ exists. In fact, consider, for
example, (\ref{eq:non_linear_system}) with
\begin{gather*}
f(x,t){=}\negthickspace\begin{bmatrix}-e^{x_{1}}+x_{2}^{3}\\
x_{1}
\end{bmatrix}\negthickspace,B^{T}(x,t){=}\negthickspace\begin{bmatrix}0\\
1
\end{bmatrix}\negthickspace,Q(x){=}\negthickspace\begin{bmatrix}-e^{x_{1}} & 1+3x_{2}^{2}\\
1+3x_{2}^{2} & 0
\end{bmatrix}\negthickspace.
\end{gather*}
 One has that the continuous function $\rho(x)=e^{-x_{1}}$ guarantees
that (\ref{eq:LMI_from_DMI}) holds for all $x\in\mathbb{R}^{2}$
and then, the system is universally exponentially stabilizable. \qed
\end{example}
One interesting consequence of considering continuity of the Finsler's
parameter with non-compact $S$  or non-polynomial $Q(s)$ and $B(s)$
 is that it opens the investigation of relaxation in more general
contexts. In fact, the property that a continuous function can be
uniformly approximated by polynomial functions on compact subsets
is what guarantees that one can reduce the search space from the set
of continuous functions to the set of polynomial functions. For more
general subsets, it is known that without compactness but under some
other technical conditions, a continuous function can be uniformly
approximated by rational functions \cite[pp.293-295]{Wil:04} or still
by functions with some other specific structure \cite{Sto:62}. 

For general sets $S$ or functions $Q\left(s\right)$ and $B\left(s\right)$,
it is very hard to check (\ref{eq:lem_sup-inf_sup-inf}). Nevertheless
bounds on $Q\left(s\right)$ and $B\left(s\right)$ can be used for
simpler tests as presented in the next lemma.
\begin{lem}
\label{lem:sup_inf_with_bounds} Let ${Q:S\rightarrow\mathbb{S}^{n}}$
and ${B:S\rightarrow\mathbb{R}^{m\times n}}$ be functions on $S\subset\mathbb{R}^{d}$
and let $\ell_{Q},\ell_{R},u_{Q},u_{R}:S\rightarrow\mathbb{R}$ such
that for all $s\in S$, 
\[
\ell_{Q}\left(s\right)I_{n}\preceq Q\left(s\right)\preceq u_{Q}\left(s\right)I_{n},\quad\ell_{R}\left(s\right)I_{n}\preceq B^{T}(s)B(s)\preceq u_{R}\left(s\right)I_{n}.
\]
Then, a necessary condition for (\ref{eq:lem_sup-inf_sup-inf}) is
\begin{equation}
\sup_{s\in S}\,\inf\left\{ \mu\in\mathbb{R}_{\geq0}\mid\ell_{Q}\left(s\right)-\mu u_{R}\left(s\right)<0\right\} <\infty\label{eq:lem_sup-inf_sup-inf_nec}
\end{equation}
and a sufficient condition for (\ref{eq:lem_sup-inf_sup-inf}) is
\begin{equation}
\sup_{s\in S}\,\inf\left\{ \mu\in\mathbb{R}_{\geq0}\mid u_{Q}\left(s\right)-\mu\ell_{R}\left(s\right)<0\right\} <\infty.\label{eq:lem_sup-inf_sup-inf_suf}
\end{equation}

In particular, if $Q$ and $B$ are scalar functions, a necessary
and sufficient condition for (\ref{eq:lem_sup-inf_sup-inf}) is

\[
Q\left(s\right)<0,\:\forall s\in S_{0}\quad and\quad\sup_{s\in S_{+}}\frac{Q(s)}{B^{2}\left(s\right)}<\infty
\]
where $S_{0}=\left\{ s\in S\mid B\left(s\right)=0\right\} $ and $S_{+}=S\setminus S_{0}.$
\end{lem}
\begin{IEEEproof}
For $\left(s,\mu\right)\in S\times\mathbb{R}$, define $\alpha\left(s,\mu\right)=\left[\ell_{Q}\left(s\right)-\mu u_{R}\left(s\right)\right]I_{n}$,
$\beta\left(s,\mu\right)=Q\left(s\right)-\mu B^{T}\left(s\right)B\left(s\right)$
and $\gamma\left(s,\mu\right)=\left[u_{Q}\left(s\right)-\mu\ell_{R}\left(s\right)\right]I_{n}$.
For each $s\in S$, we have that $\alpha\left(s,\mu\right)\preceq\beta\left(s,\mu\right)\preceq\gamma\left(s,\mu\right)$,
for all $\mu\ge0$. Therefore
\begin{equation}
\begin{aligned}\left\{ \mu\geq0\mid\gamma\left(s,\mu\right)\prec0\right\}  & \subseteq\left\{ \mu\geq0\mid\beta\left(s,\mu\right)\prec0\right\} \\
 & \subseteq\left\{ \mu\geq0\mid\alpha\left(s,\mu\right)\prec0\right\} .
\end{aligned}
\label{eq:subsets}
\end{equation}

{[}\emph{Sufficiency}{]} From (\ref{eq:subsets}), $+\infty>\sup_{s\in S}\,\inf\{\mu\geq0\mid\gamma\left(s,\mu\right)\prec0\}\geq\sup_{s\in S}\,\inf\{\mu\geq0\mid\beta\left(s,\mu\right)\prec0\}\geq\sup_{s\in S}\,\inf\{\mu\in\mathbb{R}\mid\beta\left(s,\mu\right)\prec0\}.$

{[}\emph{Necessity}{]} Since $\sup_{s\in S}\,\inf\left\{ \mu\in\mathbb{R}\mid\beta\left(s,\mu\right)\prec0\right\} <+\infty$,
exists $m>0$ such that $\inf\left\{ \mu\in\mathbb{R}\mid\beta\left(s,\mu\right)\prec0\right\} <m,\,\forall s$.
Therefore, $m\in\left\{ \mu\in\mathbb{R}\mid\beta\left(s,\mu\right)\prec0\right\} ,\,\forall s$.
Since $m>0$, we also have that $m\in\left\{ \mu\geq0\mid\beta\left(s,\mu\right)\prec0\right\} ,\,\forall s$.
From (\ref{eq:subsets}) it follows that $m\in\left\{ \mu\geq0\mid\alpha\left(s,\mu\right)\prec0\right\} ,\,\forall s$.
Therefore, $\inf\left\{ \mu\geq0\mid\alpha\left(s,\mu\right)\prec0\right\} \leq m,\,\forall s$
and $\sup_{s\in S}\,\inf\left\{ \mu\geq0\mid\alpha\left(s,\mu\right)\prec0\right\} \leq m<+\infty$.
\end{IEEEproof}
Different from the previous results, the next theorem presents a simple
case where it is possible to assure the existence of a solution. 
\begin{thm}
\label{thm:nec_suf_cond_using_bounds} Let ${Q:S\rightarrow\mathbb{S}^{n}}$
and ${B:S\rightarrow\mathbb{R}^{m\times n}}$ be matrix valued functions
on $S\subseteq\mathbb{R}^{d}$. If there are functions $\ell_{R},u_{Q}:S\rightarrow\mathbb{R}$
such that for all $s\in S$,
\begin{gather*}
0\prec\ell_{R}\left(s\right)I_{n}\preceq B^{T}(s)B(s),\quad Q\left(s\right)\preceq u_{Q}\left(s\right)I_{n},
\end{gather*}
then the statement $\FinScalarP$ holds. If further, $\ell_{R}$ and
$u_{Q}$ are continuous then the statement $\FinScalarC$ holds. Moreover,
if further $\ell_{R}$ and $u_{Q}$ are continuous and $S$ is compact,
then all the statements from $\FinScalarP$ to $\FinScalarU$ hold.
In particular, one solution to $\FinScalarU$ is
\[
\bar{\mu}=\sup_{s\in S}\frac{u_{Q}(s)+|u_{Q}(s)|}{\ell_{R}(s)}+1.
\]
\end{thm}
\begin{IEEEproof}
Since $\ell_{R}\left(s\right)>0$, it follows that $B(s)$ is full
column rank for each $s\in S$, $B^{T}(s)B(s)\in\mathbb{S}_{>0}^{n}$.
Thus, the function $\mu\left(s\right)$ obtained by extending Lemma~\ref{lem:sufficient_condition_for_mu_continuous}
pointwisely is such that ${Q(s)-\mu\left(s\right)B^{T}(s)B(s)\prec0}$
and therefore $\FinScalarP$ holds. If the functions $\ell_{R}$ and
$u_{Q}$ are continuous, then $\mu\left(s\right)$ is continuous on
$S$. If in addition $S$ is compact, then by Theorem~\ref{thm:Finsler_lemma_for_continuous_mu}
the statements from $\FinScalarP$ to $\FinScalarU$ are all equivalent. 
\end{IEEEproof}
Theorem \ref{thm:nec_suf_cond_using_bounds} can be used for determining
continuous solution $\mu\left(s\right)$ even in the case where the
matrix valued functions $Q\left(s\right)$ and $B\left(s\right)$
are not continuous. It is only necessary to find continuous bounding
functions $\ell_{R}\left(s\right)$ and $u_{Q}\left(s\right)$. In
the case where the functions $Q\left(s\right)$ and $B\left(s\right)$
are continuous, one can choose $\ell_{R}\left(s\right)=\lambda_{min}(B^{T}(s)B(s))$,
$u_{Q}\left(s\right)=\lambda_{max}(Q(s))$ and then, the following
corollary can be stated.
\begin{cor}
\label{thm:full_col_rank_constant} Let ${Q:S\rightarrow\mathbb{S}^{n}}$
and ${B:S\rightarrow\mathbb{R}^{m\times n}}$ be matrix valued functions
on $S\subseteq\mathbb{R}^{d}$ with $B(s)$ full column rank for every
$s\in S$. Then the statement $\FinScalarP$ holds. If further, $Q$
and $B$ are continuous then the statement $\FinScalarC$ holds. Moreover,
if further $Q$ and $B$ are continuous and $S$ is compact, then
all the statements from $\FinScalarP$ to $\FinScalarU$ hold. In
particular, one solution to $\FinScalarU$ is 
\[
\bar{\mu}=\sup_{s\in S}\frac{\lambda_{max}(Q(s))+|\lambda_{max}(Q(s))|}{\lambda_{min}(B^{T}(s)B(s))}+1.
\]
\end{cor}
Other combinations considering constant, polynomial, etc. functional
dependence of bounding functions $\ell_{R}\left(s\right)$ and $u_{Q}\left(s\right)$
can be trivially obtained. In the next section we comment on how the
above theorems can be applied to give simpler functional solutions
$X\left(s\right)$ for the PD-LMI in $\FinExtMatrix$.

\subsection{\label{subsec:results_for_matrix_X} Consequences for the matrix-valued
function $X\left(s\right)$}

In many control problems, one may be led to a PD-LMI in the $\FinExtMatrix$
formalism, that is, to the problem of finding a matrix function $X:S\rightarrow\mathbb{R}^{n\times m}$
satisfying $\FinExtMatrix$. 

It is known that the general results of Theorem \ref{thm:Finsler_lemma_for_matrix_polynomials}
and Theorem \ref{thm:Bliman} are valid for $\FinExtMatrix$ in the
sense that, under their hypotheses, if there is a solution to $\FinExtMatrix$
then there also exists a rational or a polynomial solution, respectively.

Remembering that if $\mu\left(s\right)$ is a solution to $\FinExtScalar$
then $X\left(s\right)=-\frac{1}{2}\mu\left(s\right)B^{T}\left(s\right)$
is a solution to $\FinExtMatrix$, it is easy to apply the results
of Section~\ref{subsec:results_for_scalar_mu} to give some sufficient
conditions that allow simple functional dependence like continuity
or polynomial dependence on $s$ for the variable $X\left(s\right)$
in $\FinExtMatrix$ without loss of generality. Among these possible
extensions, one may point out the next theorem which deals with a
case closely related to \cite{Bli:04b} and \cite{Cim:15}.
\begin{thm}
\label{thm:poly_degree_bound} Let ${Q:S\rightarrow\mathbb{S}^{n}}$
and ${B:S\rightarrow\mathbb{R}^{m\times n}}$ be matrix valued functions
on $S\subseteq\mathbb{R}^{d}$. Suppose that $B\left(s\right)$ is
a polynomial matrix.

If $S=\mathbb{R}^{d}$ and $Q$ is polynomial, it follows that if
there exists a solution $X:\mathbb{R}^{d}\rightarrow\mathbb{R}^{n\times m}$
to $\FinExtMatrix$, then there is also a solution $\bar{X}:\mathbb{R}^{d}\rightarrow\mathbb{R}^{n\times m}$
that is a polynomial matrix.

If $S\subset\mathbb{R}^{d}$ is compact, $Q$ is continuous and $B$
is polynomial of degree $g$, it follows that if there exists a solution
$X:S\rightarrow\mathbb{R}^{n\times m}$ to $\FinExtMatrix$, then
there is also a solution $\bar{X}:S\rightarrow\mathbb{R}^{n\times m}$
that is a polynomial matrix of degree $g$. In particular, if $g=0$,
then $\bar{X}$ can be taken as constant. 
\end{thm}
\begin{IEEEproof}
If there exists a matrix valued function $X:S\rightarrow\mathbb{R}^{n\times m}$
such that $\FinExtMatrix$ holds, then it follows directly from Lemma~\ref{lem:Finsler_pointwise}
that there exists a function $\mu:S\rightarrow\mathbb{R}$ satisfying
$\FinScalarP$. 

Suppose that $S=\mathbb{R}^{d}$ and $Q$ is polynomial. Since $Q$
and $B$ are polynomial matrices over $\mathbb{R}^{d}$, by Theorem~\ref{thm:Finsler_lemma_for_matrix_polynomials}
it follows that there exists a polynomial $\bar{\mu}$ satisfying
$\FinScalarPoly$. By taking $\bar{X}(s)=-\frac{1}{2}\bar{\mu}(s)B^{T}(s)$,
we have that $\bar{X}$ is a polynomial matrix. 

Suppose now that $S\subset\mathbb{R}^{d}$ is compact. Since $Q$
and $B$ are continuous functions and $S$ is compact, by Theorem~\ref{thm:cont_Q_B}
it follows that there exists a constant $\bar{\mu}\in\mathbb{R}$
satisfying $\FinScalarU$. The result follows now by taking $\bar{X}(s)=-\frac{1}{2}\bar{\mu}B^{T}(s)$.
\end{IEEEproof}
In contrast to \cite{Cim:15}, this theorem presents a case where
beside knowing the existence of a polynomial solution, it is also
possible to define its degree. This result is useful to reduce the
search space to polynomial solutions of degree less than or equal
to some specific degree, as the next example illustrates. The example
concerns how Finsler's lemma is applied along relaxation techniques
to transform a PD-LMI into a set of LMIs. This technique has been
applied, for instance in \cite{EPA:15}, where the effectiveness of
this approach is illustrated in the evaluation of LMI methods for
robust performance analysis of closed-loop longitudinal dynamics of
a civil aircraft.
\begin{example}
Consider a linear system with polytopic uncertainties, that is,
\begin{equation}
\dot{x}=A(\alpha)x(t),\label{eq:uncertain_polytopic_system}
\end{equation}
with $A(\alpha)=\sum_{i=1}^{N}\alpha_{i}A_{i}$, where $A_{i}\in\mathbb{R}^{n\times n}$
are the vertices of $A(\alpha)$ and $\alpha\in\Delta_{N}\coloneqq\left\{ \theta\in\mathbb{R}^{N}\mid\sum_{i=1}^{N}\theta_{i}=1,\theta_{i}\ge0\right\} $. 

The system~(\ref{eq:uncertain_polytopic_system}) is robustly stable
if and only if there exists $P:\Delta_{N}\rightarrow\mathbb{S}^{n}$
satisfying $P(\alpha)\succ0$ and
\begin{equation}
A^{T}(\alpha)P(\alpha)+P(\alpha)A(\alpha)\prec0\label{eq:Lyapunov_LMI-1}
\end{equation}
 for all $\alpha\in\Delta_{N}$. In order to relax this PD-LMI into
a finite set of LMIs one can impose a polynomial structure with degree
$g$ in $P(\alpha)$ \cite{OP:07a} (note that there is no loss of
generality in using a polynomial structure by Theorem~\ref{thm:Bliman}).
For instance, considering that $g=1$, $P(\alpha)=\sum_{i=1}^{N}\alpha_{i}P_{i}$
yields the PD-LMIs $\sum_{i=1}^{N}\alpha_{i}P_{i}\succ0$ and
\[
\sum_{i=1}^{N}\alpha_{i}^{2}(A_{i}^{T}P_{i}+P_{i}A_{i})+\sum_{i=1}^{N-1}\sum_{j=i+1}^{N}\alpha_{i}\alpha_{j}\left(A_{i}^{T}P_{j}+P_{j}A_{i}\right)\prec0
\]
and, consequently, a sufficient condition to the robust stability
of system~(\ref{eq:uncertain_polytopic_system}) is that the following
set of LMIs is satisfied:
\[
P_{i}\succ0,\ A_{i}^{T}P_{j}+P_{j}A_{i},\ i,j=1,\ldots,N.
\]
As shown in \cite{OP:07a}, in order to obtain a possible less conservative
set of LMIs that implies the robust stability of (\ref{eq:uncertain_polytopic_system}),
instead of applying the relaxation procedure to the standard Lyapunov
PD-LMI (\ref{eq:Lyapunov_LMI-1}), the relaxation can be applied by
imposing a polynomial structure with degree $g$ in the variables
$P(\alpha)$ and $X(\alpha)$ of the PD-LMIs $P(\alpha)\succ0$ and
\[
\begin{bmatrix}0 & P(\alpha)\\
P(\alpha) & 0
\end{bmatrix}+X(\alpha)\begin{bmatrix}A(\alpha) & -I\end{bmatrix}+\begin{bmatrix}A^{T}(\alpha)\\
-I
\end{bmatrix}X^{T}(\alpha)\prec0.
\]
In this case, there are
\[
n(5n+1)(N+g-1)!/(2g!(N-1)!)
\]
LMI scalar variables \cite{OP:07a}. However, by using Theorem~\ref{thm:poly_degree_bound}
the search of a polynomial solution to $X(\alpha)$ can be reduced,
without loss of generality, to the search of linear solutions. This
leads to
\[
n(n+1)(N+g-1)!/(2g!(N-1)!)+2n^{2}N
\]
scalar variables, which represents a search space reduction of $n^{2}O(N^{g})$
scalar variables.
\end{example}

\section{Conclusion\label{sec:conclusion}}

In this technical note, it was proposed a set of sufficient conditions
which allows the use of simple structures in the extra variable introduced
by Finsler's lemma in the context of parameter dependent matrices.
Since an unnecessarily complicated structure of parameter dependent
variables increases the computational burden without reducing the
conservatism, the results of this note may contribute to the reduction
of the computational costs for stability analysis, controller and
filter design, and any other parameter dependent results that can
be built upon the Finsler's lemma for parameter dependent systems.
While some progress has been made in the extra scalar variable, various
important issues ask for further efforts. Among them, conditions for
smooth or analytic scalar variables, more general and less conservative
conditions for continuous or polynomial solutions for the extra matrix
variable, and extended versions of the Finsler's lemma such as the
Projection Lemma, seem essential for extending optimization, control
and filter design techniques for parameter dependent systems.

\section*{Appendix}
\begin{lem}
\label{lem:sufficient_condition_for_mu_continuous} Let $M\in\mathbb{S}^{n}$
and $N\in\mathbb{S}_{>0}^{n}$. Then there exists $\mu\in\mathbb{R}$
such that $M-\mu N\prec0.$ In fact, one such $\mu$ is given by
\[
\mu=\frac{u_{M}+|u_{M}|}{\ell_{N}}+1,
\]
where $\ell_{N}$ and $u_{M}$ are any real numbers such that $0\prec\ell_{N}I_{n}\preceq N$
and $M\preceq u_{M}I_{n}$. In particular, for $\ell_{N}=\lambda_{min}(N)$
and $u_{M}=\lambda_{max}(M)$ one has
\[
\mu=\frac{\lambda_{max}(M)+|\lambda_{max}(M)|}{\lambda_{min}(N)}+1.
\]
\end{lem}
\begin{IEEEproof}
Since $N\succ0$, one has that there exists $\ell_{N}\in\mathbb{R}$
such that $0\prec\ell_{N}I_{n}\preceq N.$ Moreover, the fact that
$M\in\mathbb{S}^{n}$ yields that there exists $u_{M}\in\mathbb{R}$
such that $M\preceq u_{M}I_{n}$ and
\begin{equation}
\ell_{N}M\preceq u_{M}\ell_{N}I_{n}.\label{eq:lemma_for_continuous_mu_eq1}
\end{equation}

Since $u_{M}\in\mathbb{R}$ and $\ell_{N}>0$, there exists $\eta>0$
such that $u_{M}+\eta\ell_{N}>0.$ Indeed, it is enough to take $\eta=\frac{|u_{M}|}{\ell_{N}}+1.$
From $\eta>0$ and $\ell_{N}>0$ it follows that
\begin{equation}
u_{M}\ell_{N}<\left(u_{M}+\eta\ell_{N}\right)\ell_{N}.\label{eq:lemma_for_continuous_mu_eq3}
\end{equation}

Inequalities (\ref{eq:lemma_for_continuous_mu_eq1}) and (\ref{eq:lemma_for_continuous_mu_eq3})
yields $\ell_{N}M\prec\left(u_{M}+\eta\ell_{N}\right)\ell_{N}I_{n}\preceq\left(u_{M}+\eta\ell_{N}\right)N.$
Since $\ell_{N}\succ0$, it is enough to take
\[
\mu=\frac{u_{M}+\eta\ell_{N}}{\ell_{N}}=\frac{u_{M}+|u_{M}|}{\ell_{N}}+1.
\]
\end{IEEEproof}

\section*{Acknowledgements}

The authors would like to thank the anonymous reviewers for the excellent
suggestions and comments which led to significant improvements of
the paper. We are grateful to Prof. Jaka Cimpri\v{c} and Igor De Sant'Ana
Fontana for useful discussions. This work received support from the
Brazilian agencies CNPq and CAPES.

\balance

\bibliographystyle{ieeetr}

\end{document}